\documentclass[11pt,twoside,reqno]{amsart}

\usepackage[OT1]{fontenc}
\usepackage{type1cm}
\usepackage{lmodern}
\usepackage{microtype}

\usepackage{mathtools}
\usepackage{amssymb}
\usepackage{enumitem}
\usepackage{xcolor}
\definecolor{RED}{rgb}{1,0,0} % for AMS uppercase front matter
\usepackage{graphicx}
\usepackage{geometry}
\usepackage[bookmarks=false]{hyperref}
\hypersetup{
  colorlinks=true,
  linkcolor=black,
  anchorcolor=black,
  citecolor=black,
  filecolor=black,
  urlcolor=black
}
\usepackage[capitalize,noabbrev,nameinlink]{cleveref}

\numberwithin{equation}{section}
\allowdisplaybreaks

\theoremstyle{plain}
\newtheorem{thm}{Theorem}[section]
\newtheorem{lemma}[thm]{Lemma}
\newtheorem{cor}[thm]{Corollary}

\newtheorem{prop}[thm]{Proposition}

\theoremstyle{remark}

\newtheorem{remark}[thm]{Remark}

\renewcommand{\epsilon}{\varepsilon}

\newcommand{\F}{\mathbb{F}}
\DeclareMathOperator{\rk}{rk}

\DeclareMathOperator{\Tr}{Tr}

\DeclareMathOperator{\Span}{span}
\DeclareMathOperator{\Mat}{Mat}

\newcommand{\vu}{\mathbf{u}}
\newcommand{\vv}{\mathbf{v}}
\newcommand{\vw}{\mathbf{w}}
\newcommand{\vx}{\mathbf{x}}
\newcommand{\vy}{\mathbf{y}}

\newcommand{\M}{\mathbf M}
\title[Low-rank matrices with three prescribed values]
{Low-rank matrices with three prescribed\\off-diagonal values over $\mathbb{F}_p$}

\author{Firdavs Rakhmonov}
\address[Firdavs Rakhmonov]{School of Mathematics and Statistics,
University of St Andrews, St Andrews, KY16 9SS, Scotland}
\email{fr52@st-andrews.ac.uk}
\thanks{F. Rakhmonov was supported by a
Leverhulme Trust Research Project Grant (RPG-2023-281).}
\author{Chi Hoi Yip}
\address[Chi Hoi Yip]{Department of Mathematics, The Hong Kong University of Science and Technology, Hong Kong}
\email{machyip@ust.hk}
\date{November 1, 2026}

\begin{document}

\thispagestyle{empty}

\begin{abstract}
For every prime $p$ and every three-element set $L\subseteq\F_p$, we
prove that, for every integer $n\geq1$, there exists a symmetric
$n\times n$ matrix over $\F_p$ with zero diagonal, off-diagonal entries
in $L$, and rank $O_{p,L}(n^{1/3})$. This answers a question of Bukh.
When $0\notin L$, a standard polynomial argument shows that the exponent
$1/3$ is optimal. The proof combines an incidence matrix factorization
with a projective-plane construction controlling line-intersection
numbers modulo $p$. A separate geometric construction in characteristic
two yields better explicit constants in certain cases, and also gives
integer lifts $L'$ of the prescribed sets $L$ and symmetric $L'$-matrices
of rational rank $O_{p,L}(n^{1/3})$. When $0\notin L'$, this improves
Bukh's general quadratic lower bound, giving matrices of order
$\Omega_{p,L}(r^3)$ and rational rank at most $r$.
\end{abstract}

\maketitle
%\tableofcontents

\section{Introduction}

Let $K$ be a field and let $L\subseteq K$ be finite. Following
Bukh~\cite{bukh-ranks}, an \emph{$L$-matrix} is a square matrix over $K$
whose diagonal entries are zero and whose off-diagonal entries belong to $L$.
Matrices with prescribed entries arise naturally in linear algebraic
methods in extremal combinatorics. For $0\notin L$ and an integer $r\geq 0$, Bukh studies the extremal function
\begin{equation*}
N_K(r,L)\coloneqq
\max\bigl\{
n\geq1 :
\text{there exists an $n\times n$ $L$-matrix over $K$
of rank at most $r$}
\bigr\}.
\end{equation*} 
This function measures the largest possible order of an $L$-matrix subject to a rank constraint. Bukh investigates it using polynomial bounds and constructions from finite geometry. If $0\in L$, the possible orders are unbounded, since the zero matrix is an $L$-matrix of every order.

For a three-element set $L\subseteq K\setminus\{0\}$, a standard
polynomial argument gives a lower bound of order $n^{1/3}$ for the rank of
an $n\times n$ $L$-matrix. Bukh asks in~\cite[Problem~7]{bukh-problems}
whether this order of magnitude can always be attained over $\F_p$: for every prime $p$ and every three-element set $L\subseteq\F_p$, do there
exist $n\times n$ $L$-matrices of rank $O_{p,L}(n^{1/3})$ for every $n\geq1$? In particular, he identifies $L=\{1,4,6\}\subseteq\F_{19}$ as the simplest open case.

We completely resolve this question.

\begin{thm}
\label{thm:main}
Let $p$ be a prime. For every three-element set $L\subseteq\F_p$,
there is a constant
$C_{p,L}>0$ such that, for every integer $n\geq1$, there exists a symmetric
$L$-matrix $\M_n\in\Mat_n(\F_p)$ satisfying
\[
  \rk_{\F_p}(\M_n)\leq C_{p,L}n^{1/3}.
\]
If $0\notin L$, the exponent $1/3$ is optimal.
\end{thm}

Indeed, every $n\times n$
$L$-matrix $\M$ over any field, with $|L|=3$ and $0\notin L$, satisfies
\begin{equation}\label{eq:intro-lower}
  n\leq\binom{\rk(\M)+3}{3}.
\end{equation}
See~\cite[Proposition~1 and Lemma~10]{bukh-ranks}; we include a proof of~\eqref{eq:intro-lower} in
Appendix~\ref{sec:lower-bound}. Thus, when $0\notin L$, Theorem~\ref{thm:main} determines the minimum possible rank up to multiplicative constants depending only on $p$ and $L$. We refer to Remark~\ref{rem:general-constant} for a quantitative upper bound on $C_{p,L}$ from our proof.

The case highlighted by Bukh admits the following explicit bound.

\begin{cor}
\label{cor:F19}
Let $L=\{1,4,6\}\subseteq\F_{19}$. For every integer $n\geq1$, there exists
a symmetric $L$-matrix $\M_n\in\Mat_n(\F_{19})$ such that
\[
  \rk_{\F_{19}}(\M_n)\leq1921n^{1/3}.
\]
\end{cor}

The constant in Corollary~\ref{cor:F19} is obtained by taking $h=8$ and
$q=16$ in the construction below. We give the calculation at the end of
Section~\ref{sec:matrices}.

The characteristic-two construction giving  Corollary~\ref{cor:F19} requires suitable powers of $2$ and proves Theorem~\ref{thm:main} directly when $2$ is a primitive root modulo $p$. To remove this restriction, we use a second construction of point sets in projective planes over finite fields, controlling line-intersection numbers modulo $p$ rather than exactly. The characteristic-two construction remains useful for its smaller explicit constants and for the rational-rank result below.

\smallskip

The characteristic-two construction also yields the following strengthening over $\mathbb{Q}$. This follows from the construction itself, rather than from the finite-field rank estimate alone.

\begin{thm}
\label{thm: rank bound in Q}
Let $p$ be a prime such that $2$ is a primitive root modulo $p$, and let $L\subseteq\mathbb F_p$ be a three-element set. Then there exist a three-element set $L'\subseteq\mathbb Z$ with $L' \bmod{p} = L$ and a constant $C_{p,L}>0$ such that, for every integer \(n\ge1\), there exists a symmetric $L'$-matrix $\mathbf{M}_n\in\Mat_n(\mathbb{Z})$ satisfying
\[
\rk_{\mathbb Q}(\mathbf{M}_n)\le C_{p,L}n^{1/3}.
\]
\end{thm}

Artin's primitive root conjecture predicts that the primes satisfying the
hypothesis of Theorem~\ref{thm: rank bound in Q} have relative natural density
\[
  \prod_{\ell\ \mathrm{prime}}
     \left(1-\frac{1}{\ell(\ell-1)}\right)=0.3739558\ldots
\]
among all primes; see~\cite{moree}. Thus, the conjecture predicts that Theorem~\ref{thm: rank bound in Q} applies to approximately $37.4\%$ of all primes. A proof that infinitely many primes satisfy this primitive root hypothesis has recently been announced in~\cite{OpenAI-Artin}, while the predicted density remains conjectural. Theorem~\ref{thm: rank bound in Q} itself is unconditional for every prime satisfying the primitive root hypothesis. The method used to prove Theorem~\ref{thm:main} for arbitrary primes does not seem to remove the primitive-root hypothesis in Theorem~\ref{thm: rank bound in Q}. Indeed, it controls line-intersection numbers only modulo $p$, whereas the rational-rank argument relies on exact intersection counts.

If $0\notin L$, the sets $L'$ constructed in the proof
consist of nonzero integers and satisfy an integer linear
relation whose coefficients sum to $1$. For these sets, Bukh's general result
\cite[Theorem~3]{bukh-ranks} gives
$N_{\mathbb Q}(r,L')=\Omega(r^2)$, whereas
Theorem~\ref{thm: rank bound in Q} strengthens this to
\[
N_{\mathbb Q}(r,L')=\Omega_{p,L}(r^3).
\]
The proof is given in Section~\ref{sec: rank bound in Q}.

\subsection*{Outline of the proof}
The characteristic-two construction begins with an auxiliary
field $\F_q$, which provides a set $\mathbf{S}\subseteq\F_q^3$ of
size $D\coloneqq (h-1)q+h$, where $h$ and $q$ are powers of two. Its intersections
with proper linear subspaces of $\F_q^3$ have sizes in $\{0,1,h\}$. We prove the required
properties directly in Section~\ref{sec:geometry}.

For an integer $k\geq 1$, we set $Q\coloneqq q^k$ and regard $\mathbf{S}$ as a subset
of $\F_Q^3$. We form the incidence matrix with rows indexed by $\vx \in \F_Q^3$, columns indexed by $(\vv,t)\in \mathbf{S}\times \F_Q$, and entries recording whether $\vv\cdot\vx=t$. Its entries are taken in $\F_p$.
Multiplying this matrix by its transpose gives a square matrix with
diagonal entries $D \bmod p$, off-diagonal entries in $\{0,1,h \bmod p\}$, and rank at most $DQ$. When $2$ is a primitive root modulo $p$, we can choose $h$ and $q$
so that an affine transformation maps these three possible off-diagonal
values to $L$ and sends $D\bmod p$ to $0$. Taking principal submatrices
gives matrices of every order.

Section~\ref{sec: rank bound in Q} uses the exact intersection counts
to prove the rational-rank theorem. Finally, for general $p$,
Section~\ref{sec:general} constructs a finite projective point set with
suitable line-intersection counts modulo $p$ and uses the same
incidence-matrix argument to complete the proof of
Theorem~\ref{thm:main}.

\subsection*{Notation}
We write $\F_q$ for the finite field with $q$ elements and
$\F_q^\times\coloneqq\F_q\setminus\{0\}$. For prime $p$, we identify $\F_p$ with $\{0,1,\dots,p-1\}$ equipped with addition and multiplication modulo $p$. For $0\leq k\leq n$, we write $\mathrm{Gr}_q(k,n)$ for the Grassmannian
of $k$-dimensional $\F_q$-linear subspaces of $\F_q^n$.
We use the standard identity
\[
  |\mathrm{Gr}_q(k,n)|=\binom{n}{k}_q,
\]
where $\binom{n}{k}_q$ is the Gaussian binomial coefficient. 

If $K$ is a field and $S\subseteq K$, then 
$\Mat_{m\times n}(S)$ denotes the set of $m\times n$ matrices with entries in $S$,
and $\Mat_n(S)\coloneqq\Mat_{n\times n}(S)$. We use bold letters for vectors and
matrices. The symbol $\rk$ denotes matrix rank, with a subscript specifying the field when needed.

For a condition $E$, its indicator $\mathbf{1}_{\{E\}}$ is $1$ if $E$ holds and $0$ otherwise. For nonnegative quantities $X$ and $Y$, the notation $X=O_{p,L}(Y)$ means that $X\leq C_{p,L}Y$ for some constant $C_{p,L}>0$ depending only on $p$ and $L$. Similarly, $X=\Omega_{p,L}(Y)$ means that $X\geq c_{p,L}Y$ for some constant $c_{p,L}>0$ depending only on $p$ and $L$.

\begin{remark}
\label{zero not in L}
If $0\in L$, the zero matrix proves
Theorem~\ref{thm:main}. Taking the standard integer representatives
of the elements of $L$ as $L'$, the zero matrix also proves
Theorem~\ref{thm: rank bound in Q}.
Hence, in the proofs below we may assume that
$L\subseteq\F_p^\times$ has three distinct elements.
\end{remark}

\subsection*{Acknowledgements}
The first author thanks Boris Bukh for noting the rational-rank consequence of the characteristic-two construction and for suggesting additional references. He also thanks Jonathan Fraser for reading the original manuscript and for his helpful suggestions.

\subsection*{Use of AI tools}
The mathematical results and original proofs were developed
by the authors without AI assistance.

\section{Choice of parameters}\label{sec:parameters}

Let $L\coloneqq\{a,b,c\}\subseteq\F_p^\times$ be a three-element set,
and put $\delta\coloneqq b-a\in \F_p^{\times}$. The following lemma provides  powers of two $h\leq q$ such
that the affine map $$\F_p \to \F_p, \qquad z\longmapsto a+\delta z$$ sends $\{0,1,h \bmod{p}\}$ to $\{a,b,c\}$ and sends
$D \bmod{p}$ to $0$, where $D\coloneqq (h-1)q+h$. We will apply this map to the entries of a matrix obtained from the geometric construction. 

\begin{lemma}\label{lem:parameters}
Suppose that $\F_p^\times=\langle 2\rangle$. Then there exist integers $1\leq s\leq m$ such that, for $h \coloneqq 2^s$ and $q\coloneqq 2^m$, we have
\begin{equation}\label{eq:residues}
  h \bmod p = \frac{c-a}{b-a},
  \qquad
  q \bmod{p} = \frac{c}{b-c}.
\end{equation}
For $D\coloneqq(h-1)q+h$, these choices satisfy the following identities in $\F_p$:
\begin{equation}\label{eq:matching}
  a+\delta=b,\qquad a+h\delta=c,\qquad a+D\delta=0.
\end{equation}
\end{lemma}

\begin{proof}[Proof of Lemma \ref{lem:parameters}]
Write
\[
  \sigma\coloneqq\frac{c-a}{b-a}\in\F_p^\times\setminus\{1\},
  \qquad
  \tau\coloneqq \frac{c}{b-c}\in\F_p^\times.
\]
Since $2$ generates $\F_p^\times$, choose
$s\in\{1,\ldots,p-2\}$ and $m_0\in\{0,\ldots,p-2\}$ with
\begin{equation*}
2^s \bmod{}p= \sigma, \qquad 2^{m_0}\bmod{p}= \tau.  
\end{equation*}
Choose a nonnegative integer
$t$ large enough that $m\coloneqq m_0+t(p-1)\geq s$. Since $2^{p-1}\bmod{p}= 1$,
the numbers $h\coloneqq 2^s$ and $q\coloneqq2^m$ satisfy~\eqref{eq:residues}.

The first two identities in~\eqref{eq:matching} follow immediately from
the definitions. For the third, all arithmetic in the following display is in $\F_p$:
\begin{equation*}
a + D\delta = a +((\sigma-1)\tau + \sigma)\delta = c + \frac{c-b}{b-a} \frac{c}{b-c} (b-a)=0. \qedhere
\end{equation*}
\end{proof}

\begin{remark}
\label{dependence of parameters on p,L}
Once $p$ and $L$ are fixed, the chosen parameters $s,m,h,q,D$ depend only on $p$ and $L$.
\end{remark}

\section{The geometric construction}\label{sec:geometry}

Throughout Sections~\ref{sec:geometry}--\ref{sec: rank bound in Q}, let $q=2^m$ and $h=2^s$, where $1\leq s\leq m$. The characteristic-two construction depends
only on these two parameters.

\subsection{An anisotropic quadratic form} The absolute trace on $\F_q$ is defined by
\[
  \Tr(z)\coloneqq z+z^2+z^{2^2}+\cdots+z^{2^{m-1}} \qquad (z\in \F_q).
\]

It is an $\F_2$-linear map from $\F_q$ onto $\F_2$ satisfying
$\Tr(z^2)=\Tr(z)$ for every $z\in \F_q$. These are standard facts about traces on finite fields; see \cite{lidl}.

Choose $\theta\in\F_q$ with $\Tr(\theta)=1$, and define the quadratic form
\begin{equation}\label{eq:quadratic}
  R(x,y)\coloneqq x^2+xy+\theta y^2.
\end{equation}

\begin{lemma}
\label{lem:quadratic}
The form $R$ has the following properties: 
\begin{enumerate}[label=\arabic*)]
    \item $R$ is anisotropic: if $R(x,y)=0$, then $(x,y)=(0,0)$.
    
    \smallskip
    
    \item For every $\lambda\in\F_q^\times$, $$|\{(x,y)\in\F_q^2:R(x,y)=\lambda\}|=q+1.$$
\end{enumerate}

\end{lemma}

\begin{proof}[Proof of Lemma \ref{lem:quadratic}]

\begin{enumerate}[label=\arabic*)]
    \item Suppose that $R(x,y)=0$. If $y=0$, then $x^2=0$, so $x=0$. If $y\neq 0$, dividing by $y^2$ gives $(x/y)^2 + x/y = \theta$. Taking traces yields $$1 = \Tr(\theta) = \Tr((x/y)^2) + \Tr(x/y) = 0,$$ a contradiction.
    \smallskip
    \item Fix $\lambda \in \F_q^{\times}$ and define 
    $X_{\lambda} \coloneqq \{(x,y)\in \F_q^2: R(x,y)=\lambda\}.$
    Consider the map $$\Phi : X_\lambda \to \mathrm{Gr}_q(1,2), \qquad \Phi(x,y)\coloneqq \Span_{\F_q}\{(x,y)\}.$$ It is not difficult to verify that $\Phi$ is a bijection, and hence
    \begin{equation*}
        |\{(x,y)\in\F_q^2:R(x,y)=\lambda\}| = \binom{2}{1}_q = q+1. \qedhere
    \end{equation*}
\end{enumerate}
\end{proof}

\subsection{Intersections with affine lines}

Choose an $s$-dimensional $\F_2$-subspace $A$ of $\F_q$. Such a subspace
exists because $\F_q$ has dimension $m$ over $\F_2$ and $1\leq s\leq m$. In particular, $|A|=h$.
Put
\begin{equation}\label{eq:affine-set}
  R^{-1}(A)\coloneqq\{(x,y)\in\F_q^2:R(x,y)\in A\}.
\end{equation}
The next lemma gives the key intersection property.

\begin{lemma}\label{lem:lines}
For every affine line $\ell\subseteq\F_q^2$,
\[
  |\ell\cap R^{-1}(A)|\in\{0,h\}.
\]
\end{lemma}

\begin{proof}[Proof of Lemma~\ref{lem:lines}]
We distinguish two cases.
\begin{enumerate}[label=\arabic*)]
    \item Suppose that $\mathbf{0}\in \ell$. Write $\ell = \{t\mathbf{v}:t \in \F_q\}$, where $\vv\neq \mathbf{0}$. By anisotropy, $R(\mathbf{v})\neq 0$. Since the map $t\mapsto t^2R(\vv)$ is a bijection on $\F_q$,
    \begin{equation*}
        |\ell \cap R^{-1}(A)| = |\{t\in \F_q: t^2R(\mathbf{v})\in A\}| = h.
    \end{equation*}
    \item Suppose that $\mathbf{0} \notin \ell$. Write $\ell = \{\mathbf{u} + t\mathbf{v}: t\in \F_q\}$, where $\mathbf{u}$ and $\mathbf{v}$ are linearly independent over $\F_q$. For every $t\in \F_q$, 
    \begin{equation*}
        R(\mathbf{u} + t\mathbf{v}) = \alpha t^2 + \beta t + \gamma,
    \end{equation*}
    where
    \begin{equation*}
        \alpha \coloneqq R(\vv) \in \F_q^{\times}, \qquad \beta\coloneqq 
u_1v_2 +u_2v_1\in \F_q^{\times}, \qquad \gamma \coloneqq R(\vu)\in \F_q^{\times}.
\end{equation*}

\smallskip

Consider the $\F_2$-linear map $$T: \F_q \to \F_q, \qquad T(t)\coloneqq \alpha t^2 + \beta t.$$
Its kernel is $\{0, \beta/\alpha\}$, so its image $V\coloneqq \operatorname{im}(T)$ is an $\F_2$-subspace of index two in $\F_q$, and every element of $V$ has exactly two preimages.

We have $\gamma \notin V$. Otherwise, $\gamma = T(t_0)$ for some $t_0\in \F_q$, and then $R(\vu + t_0\vv)=T(t_0)+\gamma = 0$. By anisotropy, $\vu + t_0\vv = \mathbf{0}$, contradicting $\mathbf{0}\notin \ell$. Consequently,
\begin{equation*}
    |\ell \cap R^{-1}(A)| = 2|(A+\gamma)\cap V|.
\end{equation*}

If $A\subseteq V$, then $(A+\gamma)\cap V = \varnothing$, because $\gamma \notin V$. Hence, $$|\ell \cap R^{-1}(A)| = 0.$$

If $A\not\subseteq V$, choose $a_0\in A\setminus V$. Then $\F_q/V = \{V, a_0+V\}$ and
\begin{equation*}
    A = (A\cap V) \sqcup (A\cap (a_0+V)).
\end{equation*}
Translation by $a_0$ is a bijection between the two parts, so $|A\cap V| = h/2$. Since $a_0$ and $\gamma$ lie in the same nonzero coset of $V$,
\begin{equation*}
|(A+\gamma)\cap V| = |A\cap V| = \frac{h}{2}.
\end{equation*}
It follows that $|\ell \cap R^{-1}(A)|=h$.
\end{enumerate}
Thus, in all cases, $|\ell \cap R^{-1}(A)| \in \{0,h\}$.
\end{proof}

\subsection{The affine lift}

Lift $R^{-1}(A)$ to the plane $z=1$ in $\F_q^3$ by defining
\begin{equation}\label{eq:lift}
  \mathbf{S}\coloneqq\{(x,y,1)\in\F_q^3: R(x,y)\in A\}.
\end{equation}

\begin{prop}\label{prop:subspaces}
The set $\mathbf{S}$ has size $D=(h-1)q+h$. Moreover, the following hold
for every proper $\F_q$-linear subspace $W$ of $\F_q^3$:
\begin{enumerate}[label=\arabic*)]
  \item If $W\in\mathrm{Gr}_q(0,3)$, then $|\mathbf{S}\cap W|=0$;
  \item If $W\in\mathrm{Gr}_q(1,3)$, then $|\mathbf{S}\cap W|\leq1$;
  \item If $W\in\mathrm{Gr}_q(2,3)$, then $|\mathbf{S}\cap W|\in\{0,h\}$.
\end{enumerate}
\end{prop}

\begin{proof}[Proof of Proposition \ref{prop:subspaces}]
By Lemma~\ref{lem:quadratic},
\[
  |\mathbf{S}|=1+(|A|-1)(q+1)=(h-1)q+h.
\]

The zero vector does not belong to $\mathbf{S}$, proving ~1). A one-dimensional subspace of $\F_q^3$ contains at most one vector
with third coordinate $1$, proving~2).

For~3), write
\[
  W=\{(x,y,z)\in\F_q^3:\alpha x+\beta y+\gamma z=0\}
\]
for some $(\alpha,\beta,\gamma)\in \F_q^3\setminus \{\mathbf{0}\}$. If
$\alpha=\beta=0$, then $W=\{z=0\}$ and $W\cap\mathbf{S}=\varnothing$.
Otherwise, $\alpha x+\beta y+\gamma=0$ defines an affine line $\ell$ in
$\F_q^2$. Therefore, $|\mathbf{S}\cap W| = |\ell\cap R^{-1}(A)|$, and Lemma~\ref{lem:lines} applies.
\end{proof}

\section{The incidence matrix and the rank bound in $\F_p$}
\label{sec:matrices}

Fix the parameters from Lemma~\ref{lem:parameters} and the set $\mathbf{S}$
from Section~\ref{sec:geometry}. For an integer $k\geq1$, put $Q\coloneqq q^k$
and identify $\F_q$ with the subfield of $\F_Q$ having $q$ elements.
The set $\mathbf{S}$ remains fixed; only the field containing its coordinates
is enlarged. For $\vv\coloneqq (v_1,v_2,v_3)$ and $\vw\coloneqq (w_1,w_2,w_3)\in\F_Q^3$, write
$$\vv\cdot\vw \coloneqq v_1w_1+v_2w_2+v_3w_3\in \F_Q.$$

\subsection{Orthogonality counts over extension fields}

For $\vw\in\F_Q^3$, define 
\[
  N(\vw)\coloneqq|\{\vv\in\mathbf{S}:\vv\cdot\vw=0\}|.
\]

\begin{lemma}\label{lem:orthogonality}
We have $N(\mathbf{0})=D$, and
\[
  N(\vw)\in\{0,1,h\}\qquad  (\vw\in\F_Q^3\setminus\{\mathbf{0}\}).
\]
\end{lemma}

\begin{proof}[Proof of Lemma \ref{lem:orthogonality}]
The first assertion follows from $|\mathbf{S}|=D$. For $\vw\ne\mathbf{0}$,
consider the $\F_q$-linear map
\[
  T_{\vw}:\F_q^3\longrightarrow\F_Q,\qquad
  T_{\vw}(\vv)\coloneqq \vv\cdot\vw.
\]
Choose a coordinate $j$ with $w_j\ne0$. The standard basis vector
$\mathbf e_j\in\F_q^3$ satisfies $T_{\vw}(\mathbf e_j)=w_j\ne0$,
so $T_{\vw}$ is nonzero. The rank--nullity theorem gives
\[
  \dim_{\F_q}\ker(T_{\vw})
   =3-\dim_{\F_q}\operatorname{im}(T_{\vw})\leq2.
\]
Since $N(\vw)=|\mathbf{S}\cap\ker(T_{\vw})|$, the result follows from
Proposition~\ref{prop:subspaces}.
\end{proof}

\subsection{The incidence factorization}\label{sec:incidence-factorization}

Define the incidence matrix with row index set $\F_Q^3$ and column index set $\mathbf{S}\times \F_Q$ by
\[
  \mathbf{B}\in\Mat_{Q^3\times DQ}(\F_p),\qquad
  \mathbf{B}_{\vx,(\vv,t)}\coloneqq \mathbf1_{\{\vv\cdot\vx=t\}}.
\]
The dot product is evaluated in $\F_Q$, whereas the
indicator is regarded as $0$ or $1$ in $\F_p$.

\begin{lemma}\label{lem:gram}
For all $\vx,\vy\in\F_Q^3$,
\begin{equation}\label{eq:gram}
  (\mathbf{B}\mathbf{B}^{\mathsf T})_{\vx,\vy}= N(\vx-\vy) \bmod{p}.
\end{equation}
\end{lemma}

\begin{proof}[Proof of Lemma \ref{lem:gram}]
By matrix multiplication,
\begin{equation*}
  (\mathbf{B}\mathbf{B}^{\mathsf T})_{\vx,\vy}
  =\sum_{\vv\in\mathbf{S}}\sum_{t\in\F_Q}
      \mathbf1_{\{\vv\cdot\vx=t\}}\,
      \mathbf1_{\{\vv\cdot\vy=t\}}=\sum_{\vv\in\mathbf{S}}\mathbf1_{\{\vv\cdot(\vx-\vy)=0\}}.\\
\end{equation*}
The products and sums are taken in \(\mathbb F_p\), so the final sum equals $N(\mathbf x-\mathbf y)\bmod p$.
\end{proof}

We first prove Theorem~\ref{thm:main} under the primitive-root hypothesis.

\begin{proof}[Proof of Theorem~\ref{thm:main} when $2$ is a primitive root modulo $p$]
By Remark~\ref{zero not in L}, it suffices to consider $0\notin L$. For an integer $k\geq 1$, put $Q\coloneqq q^k$, let $\mathbf{J}$ be the all-ones
matrix of order $Q^3$, and set
\begin{equation}\label{eq:matrix}
\M\coloneqq a\mathbf{J}+\delta\mathbf{B}\mathbf{B}^{\mathsf T}\in\Mat_{Q^3}(\F_p).
\end{equation}
This matrix is symmetric. By Lemmas~\ref{lem:orthogonality}
and~\ref{lem:gram}, its diagonal entries are
\[
  \M_{\vx,\vx}=a+D\delta=0,
\]
where the last equality follows from~\eqref{eq:matching}. For $\vx\ne\vy$,
the same lemmas and~\eqref{eq:matching} give
\[
  \M_{\vx,\vy}
    =a+ N(\vx-\vy)\delta \in \{a,a+\delta,a+h\delta\}=\{a,b,c\}=L.
\]
Thus, $\M$ is a symmetric $L$-matrix. Since $\mathbf{B}$ has $DQ$ columns,
\begin{equation}\label{eq:rank}
  \rk_{\F_p}(\M)
   \leq\rk_{\F_p}(a\mathbf{J})+\rk_{\F_p}(\delta\mathbf{B}\mathbf{B}^{\mathsf T})
   \leq1+\rk_{\F_p}(\mathbf{B})
   \leq1+DQ.
\end{equation}
To obtain a matrix of arbitrary order $n\geq1$, choose the smallest
integer $k\geq1$ with $q^{3k}\geq n$, and put $Q=q^k$.
If $k=1$, then $Q=q\leq qn^{1/3}$. If $k>1$, minimality gives
$q^{3(k-1)}<n$, and hence $Q<qn^{1/3}$. Thus, in either case,
\[
  Q\leq qn^{1/3}.
\]
Take any principal $n\times n$ submatrix $\M_n$ of $\M$. It remains
a symmetric $L$-matrix, and its rank is at most that of $\mathbf{M}$. Hence ~\eqref{eq:rank} yields
\begin{equation}\label{eq:explicit-bound}
  \rk_{\F_p}(\M_n)
   \leq1+DQ
   \leq1+Dq\,n^{1/3}
   \leq(1+Dq)n^{1/3}.
\end{equation}
The constant $1+Dq$ depends only on $p$ and $L$; see Remark \ref{dependence of parameters on p,L}. The optimality of the
exponent $1/3$ follows from~\eqref{eq:intro-lower}, proved in
Appendix~\ref{sec:lower-bound}.
\end{proof}

\begin{proof}[Proof of Corollary~\ref{cor:F19}]
Take $a=1$, $b=4$, and $c=6$ in $\F_{19}$. Then
$\delta=b-a=3$, and
\[
  \frac{c-a}{b-a}=8=2^3,
  \qquad
  \frac{c}{b-c}=16=2^4
  \quad\text{in }\F_{19}.
\]
We may therefore use $h=8$ and $q=16$. These choices give
$D=(8-1)\times 16+8=120$. The construction above
and~\eqref{eq:explicit-bound} give the constant
$1+Dq=1+120\times 16=1921$.
\end{proof}

\section{Proof of Theorem~\ref{thm: rank bound in Q}}
\label{sec: rank bound in Q}

Suppose that $\F_p^{\times} = \langle 2\rangle$ and $L=\{a,b,c\}\subseteq \F_p^{\times}$. Lemma~\ref{lem:parameters} provides integer parameters $h,q,D$ satisfying \eqref{eq:residues} and \eqref{eq:matching}. Recall that $\delta\coloneqq b-a\in\F_p^\times$, and choose a
nonzero integer representative of $\delta$.
All arithmetic involving the matrix entries below is performed
in $\mathbb{Z}$.

Define the integers
\begin{equation*}
    \alpha\coloneqq -D\delta ,\qquad \beta \coloneqq \delta-D\delta, \qquad \gamma\coloneqq h\delta - D\delta.
\end{equation*}
Since $D>h>1$ and $\delta \neq 0$, these are nonzero and pairwise distinct. Set
\begin{equation*}
    L'\coloneqq \{\alpha, \beta, \gamma\}.
\end{equation*}
Reducing $\alpha$, $\beta$, and $\gamma$ modulo $p$ and
using~\eqref{eq:matching}, we obtain
\begin{equation*}
    \alpha \bmod{p} = a,\qquad \beta \bmod{p}=b,\qquad \gamma \bmod{p} = c.
\end{equation*}
Thus, $L'\bmod{p} = L$. Moreover,
$(1-D)\alpha+D\beta=0$ and $(1-D)+D=1$, so the elements of
$L'$ satisfy the integer linear relation described in the introduction.

For an integer $k\geq 1$, put $Q=q^k$ and use the same incidence construction as in Section~\ref{sec:matrices}, now with entries in $\mathbb{Z}$:
\[
  \mathbf{B}\in\Mat_{Q^3\times DQ}(\mathbb{Z}),\qquad
  \mathbf{B}_{\vx,(\vv,t)}\coloneqq \mathbf1_{\{\vv\cdot\vx=t\}}.
\]
Here the rows are indexed by $\vx\in \F_Q^3$, the columns by $(\vv,t)\in \mathbf{S}\times \F_Q$, and the dot product is evaluated in $\F_Q$. The indicator is regarded as $0$ or $1$ in $\mathbb{Z}$. The counting argument from Lemma~\ref{lem:gram}, performed over $\mathbb{Z}$, gives
\begin{equation*}
    (\mathbf{B}\mathbf{B}^{\mathsf T})_{\vx,\vy} = N(\vx-\vy) \qquad (\vx,\vy \in \F_Q^3).
\end{equation*}
Let $\mathbf{J}$ be the all-ones matrix of order $Q^3$, and set 
\begin{equation*}
\M\coloneqq -D\delta\mathbf{J}+\delta\mathbf{B}\mathbf{B}^{\mathsf T}\in\Mat_{Q^3}(\mathbb{Z}).
\end{equation*}
By Lemma~\ref{lem:orthogonality},
\begin{equation*}
    \mathbf{M}_{\vx,\vx} = 0, \qquad \mathbf{M}_{\vx,\vy} \in \{\alpha,\beta,\gamma\}=L' \qquad (\vx\neq \vy). 
\end{equation*}
Hence $\mathbf{M}$ is a symmetric $L'$-matrix, and the same rank estimate as in Section~\ref{sec:matrices} gives
\begin{equation*}
    \rk_{\mathbb{Q}}(\mathbf{M})\leq 1 + DQ.
\end{equation*}
For any $n\geq 1$, choose $k$ and take a principal $n\times n$ submatrix as in the proof of Theorem~\ref{thm:main}. This gives a symmetric $L'$-matrix $\mathbf{M}_n\in \text{Mat}_n(\mathbb{Z})$ satisfying
\begin{equation*}
    \rk_{\mathbb{Q}}(\mathbf{M}_n) \leq (1 + Dq)n^{1/3}.
\end{equation*}
Since $D$ and $q$ depend only on $p$ and $L$, taking $C_{p,L}\coloneqq 1 + Dq$ proves Theorem~\ref{thm: rank bound in Q}.

Taking $n=\lfloor(r/C_{p,L})^3\rfloor$ for sufficiently large \(r\) (say $r\geq 2C_{p,L}$), we obtain
\[
N_{\mathbb Q}(r,L')\geq n = \Omega_{p,L}(r^3).
\]
Thus, for the sets \(L'\) constructed above, our result improves the general quadratic lower bound supplied by \cite[Theorem~3]{bukh-ranks} to a cubic one.

\section{Proof of Theorem~\ref{thm:main}}\label{sec:general}

We now complete the proof of Theorem~\ref{thm:main} for arbitrary primes.
By Remark~\ref{zero not in L}, we may assume that
$L=\{a,b,c\}\subseteq\F_p^\times$, so $p\geq5$. Put
\begin{equation}\label{eq:general-parameters}
\delta=b-a,\qquad
\eta=\frac{c-a}{b-a},\qquad
 d_0=-\frac{a}{b-a}\quad\text{in }\F_p.
\end{equation}
Then $\eta\notin\{0,1\}$ and $d_0\notin\{0,1,\eta\}$. The affine map
$z\mapsto a+\delta z$ sends $0,1,\eta$ to $a,b,c$, respectively,
and sends $d_0$ to $0$.

In what follows, $r$ denotes a prime power. We write $\mathbb{P}^2(\F_r)$ for the projective plane over $\F_r$. We use homogeneous coordinates for points in $\mathbb{P}^2(\F_r)$. An $\F_r$-line in $\mathbb{P}^2(\F_r)$ is a projective line defined over $\F_r$.
The characteristic of $\F_r$ need not be $p$, while all congruences involving point counts are taken modulo $p$.

The incidence construction in Section~\ref{sec:matrices} has the following
projective formulation, which we record for use below.

\begin{lemma}
\label{lem:general-incidence}
Suppose that $S\subseteq\mathbb P^2(\F_r)$ satisfies
\[
 |S|\equiv d_0\pmod p,\qquad
 |S\cap\ell|\in\{0,1,\eta\}\pmod p
 \quad\text{for every $\F_r$-line $\ell$.}
\]
Then, for every integer $n\geq1$, there is a symmetric $n\times n$
$L$-matrix over $\F_p$ of rank at most $(1+r|S|)n^{1/3}$.
\end{lemma}

\begin{proof}[Proof of Lemma~\ref{lem:general-incidence}]
For each $P\in S$, choose homogeneous coordinates
$P=[a_P:b_P:c_P]$ over $\F_r$. Fix an integer $j\geq1$,
put $Q=r^j$, and define a $0$--$1$ matrix $\mathbf B$ over $\F_p$.
Its rows are indexed by $\vx=(x_1,x_2,x_3)\in\F_Q^3$ and its
columns by $(P,u)\in S\times\F_Q$, with
\[
 \mathbf B_{\vx,(P,u)}
  =\mathbf1_{\{a_Px_1+b_Px_2+c_Px_3=u\}}.
\]
For $\vx,\vy\in\F_Q^3$, matrix multiplication gives
\begin{equation}\label{eq:general-gram}
 (\mathbf B\mathbf B^{\mathsf T})_{\vx,\vy}
 =\bigl|\{P\in S:
    a_P(x_1-y_1)+b_P(x_2-y_2)+c_P(x_3-y_3)=0\}\bigr|\pmod p.
\end{equation}
When $\vx=\vy$, this is $|S|\equiv d_0\pmod p$.
When $\vx\ne\vy$, the equation in~\eqref{eq:general-gram}
defines a projective line $\ell$ over $\F_Q$. If $\ell$ contains
two distinct $\F_r$-points, it is the extension of the
$\F_r$-line through those points, and all its $\F_r$-points lie on
that line. Otherwise, $\ell$ contains at most one point of $S$.
It follows that every off-diagonal entry of
$\mathbf B\mathbf B^{\mathsf T}$ belongs to $\{0,1,\eta\}$ in $\F_p$.

Let $\mathbf J$ be the all-ones matrix of order $Q^3$ and set
\[
 \M=a\mathbf J+\delta\mathbf B\mathbf B^{\mathsf T}.
\]
By~\eqref{eq:general-parameters}, its diagonal entries are
$a+\delta d_0=0$, while its off-diagonal entries belong to
$\{a,a+\delta,a+\delta\eta\}=L$.
Thus $\M$ is a symmetric $L$-matrix. As in Section~\ref{sec:incidence-factorization}, the $|S|Q$
columns of $\mathbf B$ give $\rk_{\F_p}(\M)\leq1+|S|Q$.
Taking an $n\times n$ principal submatrix with $j$ minimal
such that $r^{3j}\geq n$ gives the claimed bound.
\end{proof}

It remains to construct a point set satisfying the hypotheses
of Lemma~\ref{lem:general-incidence}. The next lemma makes this
adjustment possible. It replaces each point
by several points on a line through a common point $V$ over an extension
field, and may also include $V$. We choose the extension degree $2p-1$
both for linear independence and to preserve the residue of the field
size modulo $p$ in the final proof.

\begin{lemma}
\label{lem:general-extension}
Let $\eta\in\F_p\setminus\{0,1\}$, and suppose that
$S\subseteq\mathbb P^2(\F_r)$ satisfies
\[
 |S\cap\ell|\in\{0,1,\eta\}\pmod p
 \quad\text{for every $\F_r$-line $\ell$.}
\]
Suppose $0\leq k\leq p-2$ and $\varepsilon\in\{0,1\}$ satisfy
\[
 r^k+\varepsilon\in\{0,1,\eta\}\pmod p.
\]
Then there exists $S'\subseteq\mathbb P^2(\F_{r^{2p-1}})$ such that
\[
 |S'|=r^k|S|+\varepsilon,\qquad
 |S'\cap\ell|\in\{0,1,\eta\}\pmod p
 \quad\text{for every $\F_{r^{2p-1}}$-line $\ell$.}
\]
\end{lemma}

\begin{proof}[Proof of Lemma~\ref{lem:general-extension}]
Put $m=2p-1$ and choose $\alpha\in\F_{r^m}$ of degree $m$
over $\F_r$. Set
\[
 U=\Span_{\F_r}\{1,\alpha,\ldots,\alpha^{k-1}\},\qquad
 V=[\alpha:\alpha^2:\alpha^3]\in\mathbb P^2(\F_{r^m}),
\]
where $U=\{0\}$ if $k=0$. Then $|U|=r^k$. Since
$1,\alpha,\ldots,\alpha^{m-1}$ are linearly independent over
$\F_r$, the point $V$ lies on no $\F_r$-line. Indeed, such a line would give a nontrivial $\F_r$-linear relation among $\alpha,\alpha^2,\alpha^3$. Moreover,
\begin{equation}\label{eq:general-avoid}
 1\notin\alpha U+\alpha^2U+\alpha^3U.
\end{equation}
Indeed, the sum on the right is contained in
$\Span_{\F_r}\{\alpha,\ldots,\alpha^{k+2}\}$, which does not
contain $1$ because $k+2\leq p<m$.

For each $P=[a_P:b_P:c_P]\in S$, choose a representative over
$\F_r$ and define
\[
 C_P=\bigl\{[a_P-u\alpha:b_P-u\alpha^2:c_P-u\alpha^3]
     :u\in U\bigr\}.
\]
The set $C_P$ is independent of the chosen representative: rescaling $P$ by an element of $\F_r^\times$ simply rescales $u$ within the $\F_r$-subspace $U$.
The coordinate vectors representing $P$ and $V$ are linearly
independent over $\F_{r^m}$. Thus $C_P$ consists of $r^k$ distinct
points, none equal to $V$, all on the line joining $P$ to $V$.
The lines joining $V$ to distinct points of $S$ are distinct,
since otherwise $V$ would lie on the $\F_r$-line joining two
of those points. Hence the sets $C_P$ are pairwise disjoint.
Set
\[
 S'=\bigcup_{P\in S}C_P\quad(\varepsilon=0),\qquad
 S'=\{V\}\cup\bigcup_{P\in S}C_P\quad(\varepsilon=1).
\]
It follows that $|S'|=r^k|S|+\varepsilon$.

We verify the line-intersection condition. An $\F_{r^m}$-line
through $V$ contains either one entire $C_P$ or none of them.
Its intersection with $S'$ therefore has size $r^k+\varepsilon$
or $\varepsilon$, both in $\{0,1,\eta\}$ modulo $p$.

Now let $\ell$ be an $\F_{r^m}$-line not containing $V$. Write
\[
 \ell:\ z_1X+z_2Y+z_3Z=0,
\]
scaling the coefficients so that
$z_1\alpha+z_2\alpha^2+z_3\alpha^3=1$. Then
\[
 [a_P-u\alpha:b_P-u\alpha^2:c_P-u\alpha^3]\in\ell
 \quad\Longleftrightarrow\quad
 u=z_1a_P+z_2b_P+z_3c_P.
\]
Thus $C_P$ meets $\ell$ in exactly one point if the expression
on the right belongs to $U$, and otherwise it is disjoint from $\ell$.
Define the $\F_r$-subspace
\[
 W=\{(x,y,z)\in\F_r^3:z_1x+z_2y+z_3z\in U\}.
\]
Because $U$ is an $\F_r$-subspace, membership in $W$ is invariant
under multiplication by nonzero $\F_r$-scalars. Hence, writing
$\mathbb P(W)$ for its projectivization, we have
\begin{equation}\label{eq:general-intersection}
 |S'\cap\ell|=|S\cap\mathbb P(W)|.
\end{equation}
We claim that $W\ne\F_r^3$. Otherwise all three standard basis
vectors belong to $W$, and hence $z_1,z_2,z_3\in U$. The
normalization of $\ell$ would then give
\[
 1=z_1\alpha+z_2\alpha^2+z_3\alpha^3
 \in\alpha U+\alpha^2U+\alpha^3U,
\]
contrary to~\eqref{eq:general-avoid}. Thus $\mathbb P(W)$ is
empty, a single $\F_r$-point, or an $\F_r$-line. In the first
two cases the right-hand side of~\eqref{eq:general-intersection}
is $0$ or $1$; in the last case it lies in $\{0,1,\eta\}$
modulo $p$ by assumption. This completes the proof.
\end{proof}

\begin{proof}[Proof of Theorem~\ref{thm:main} in general]
Choose a generator $\gamma$ of $\F_p^\times$. By Dirichlet's
theorem, there is a prime $r_0\equiv\gamma\pmod p$.
Let $\mathcal R\subseteq\F_p$ consist of the residues $|S|\bmod p$
attainable by sets $S\subseteq\mathbb P^2(\F_r)$ for some prime
power $r$ satisfying
\[
 r\equiv\gamma\pmod p,\qquad
 |S\cap\ell|\in\{0,1,\eta\}\pmod p
 \quad\text{for every $\F_r$-line $\ell$.}
\]
The empty set in $\mathbb P^2(\F_{r_0})$ shows that
$0\in\mathcal R$.

Choose $k_0,k_1\in\{0,\ldots,p-2\}$ with
$\gamma^{k_0}=\eta$ and $\gamma^{k_1}=\eta-1$, and define
\[
 f(x)=\eta x,\qquad g(x)=(\eta-1)x+1.
\]
We claim that $f(\mathcal R),g(\mathcal R)\subseteq\mathcal R$.
Indeed, if $d\in\mathcal R$ is represented by $S$ over $\F_r$,
then $r^{k_0}\equiv\eta\pmod p$. Lemma~\ref{lem:general-extension}
with $(k,\varepsilon)=(k_0,0)$ gives a set $S'$ whose cardinality
is congruent to
\[
 |S'|=r^{k_0}|S|\equiv\eta d=f(d)\pmod p.
\]
The permitted line-intersection residues are preserved, and
\[
 r^{2p-1}\equiv\gamma^{2p-1}=\gamma\pmod p,
\]
so $S'$ also satisfies the condition on the field size. Hence
$f(d)\in\mathcal R$. Similarly, using $(k,\varepsilon)=(k_1,1)$
produces a set $S''$ with
\[
 |S''|=r^{k_1}|S|+1
 \equiv(\eta-1)d+1=g(d)\pmod p,
\]
and again the field size is congruent to $\gamma$ modulo $p$.
Thus $g(d)\in\mathcal R$ as well.

Since $\eta$ and $\eta-1$ are nonzero, $f$ and $g$ are
permutations of $\F_p$. As they preserve the finite nonempty set
$\mathcal R$, they also permute $\mathcal R$. Put
$\sigma=\sum_{d\in\mathcal R}d\in\F_p$. Summing the images of
$\mathcal R$ under the two permutations gives
\[
 \sigma=\eta\sigma,\qquad
 \sigma=(\eta-1)\sigma+|\mathcal R|.
\]
The first identity yields $\sigma=0$, since $\eta\ne1$.
The second then implies $|\mathcal R|\equiv0\pmod p$.
Since $1\leq|\mathcal R|\leq p$, it follows that
$\mathcal R=\F_p$.

In particular, there exist a prime power $r$ and a point set
$S\subseteq\mathbb P^2(\F_r)$ satisfying
\[
 |S|\equiv d_0\pmod p,\qquad
 |S\cap\ell|\in\{0,1,\eta\}\pmod p
 \quad\text{for every $\F_r$-line $\ell$.}
\]
Applying Lemma~\ref{lem:general-incidence} to this point set
gives, for every $n\geq1$, a symmetric $n\times n$ $L$-matrix
of rank at most $(1+r|S|)n^{1/3}$. This completes the proof.
\end{proof}

\begin{remark}\label{rem:general-constant}
When $0\notin L$, the construction also gives an explicit, though
very large, bound for $C_{p,L}$. Let $g^*(p)$ be the least prime
that is a primitive root modulo $p$. In the preceding proof,
choose $\gamma\equiv g^*(p)\pmod p$ and take $r_0=g^*(p)$.
The bound $g^*(p)\ll p^{3.1}$ of Ha~\cite{ha} gives a uniform
polynomial estimate for the starting field size.

Let $\mathcal T$ be the set of residues reachable from $0$ by
finitely many applications of $f$ and $g$, including the empty
sequence. Then $0\in\mathcal T$ and
$f(\mathcal T),g(\mathcal T)\subseteq\mathcal T$. Since the maps
are injective, they permute the finite set $\mathcal T$.
Writing $\tau=\sum_{d\in\mathcal T}d$, we obtain
\[
 \tau=\eta\tau,\qquad
 \tau=(\eta-1)\tau+|\mathcal T|.
\]
As in the proof above, these identities imply $\mathcal T=\F_p$.
In particular, there is a sequence
\[
 0=x_0,x_1,\ldots,x_t=d_0,\qquad
 x_{i+1}\in\{f(x_i),g(x_i)\}.
\]
Choose one of minimal length $t$. Its residues are distinct,
since otherwise a repeated segment could be deleted, so
$t\leq p-1$. Starting with the empty point set over
$\F_{r_0}$, realize each step by
Lemma~\ref{lem:general-extension}. Because the field size
remains congruent to $\gamma$ modulo $p$ at every step, the
resulting set lies in $\mathbb P^2(\F_r)$ for some
\[
 r=r_0^{(2p-1)^t}
 \leq\bigl(g^*(p)\bigr)^{(2p-1)^{p-1}}.
\]
Since $|S|\leq r^2+r+1$, the rank estimate at the end of the
proof of Theorem~\ref{thm:main} gives
\begin{equation}\label{eq:general-constant}
 C_{p,L}\leq2\bigl(g^*(p)\bigr)^{3(2p-1)^{p-1}}
 \leq\exp \bigl(O((2p-1)^{p-1}\log p)\bigr),
\end{equation}
with an absolute implied constant.

Under the generalized Riemann Hypothesis (GRH), the estimates of Shoup~\cite{shoup} and Martin~\cite{martin} yield
\[
 g^*(p)\ll
 \bigl(\omega(p-1)\log(2\omega(p-1))\bigr)^4(\log p)^2
 \ll(\log p)^6,
\]
where $\omega(m)$ is the number of distinct prime divisors
of $m$. Consequently,~\eqref{eq:general-constant} improves to
\[
 C_{p,L}\leq
 \exp\!\bigl(O((2p-1)^{p-1}\log\log p)\bigr).
\]
The characteristic-two construction in
Sections~\ref{sec:parameters}--\ref{sec:matrices} can give
much smaller constants for particular choices of $p$ and $L$;
it also provides the exact intersection counts used in the
rational-rank theorem of Section~\ref{sec: rank bound in Q}.
\end{remark}

\appendix
\section{The matching lower bound}\label{sec:lower-bound}

We record the standard polynomial argument behind~\eqref{eq:intro-lower}.
It works over an arbitrary field.

\begin{prop}
\label{prop:lower}
Let $K$ be a field and let $L=\{a,b,c\}\subseteq K\setminus\{0\}$
have three distinct elements. Every $L$-matrix $\M\in\Mat_n(K)$
satisfies
\[
  n\leq\binom{\rk_K(\M)+3}{3}.
\]
In particular, $\rk_K(\M)\geq(6n)^{1/3}-3$.
\end{prop}

\begin{proof}[Proof of Proposition~\ref{prop:lower}]
Put $r=\rk_K(\M)$. If $r=0$, then $\M$ is the zero matrix, and
$0\notin L$ forces $n=1$. Both conclusions hold in this case.
Suppose that $r\geq1$ and choose a rank factorization
\[
  \M=\mathbf U\mathbf V,\qquad
  \mathbf U\in\Mat_{n\times r}(K),\qquad
  \mathbf V\in\Mat_{r\times n}(K).
\]
Let $\vu_i$ be the $i$th row of $\mathbf U$ and $\vv_j$ the $j$th
column of $\mathbf V$. Then $\mathbf{M}_{ij}=\vu_i\cdot\vv_j$.

Define $P(x)\coloneqq (x-a)(x-b)(x-c)\in K[x]$ and observe that $P(0)=-abc \neq 0$. For $1\leq i\leq n$, define
\[
  f_i(X_1,\ldots,X_r)
   \coloneqq P(u_{i1}X_1+\cdots+u_{ir}X_r).
\]
Each $f_i$ is a polynomial of total degree at most three. Since $\M$
is an $L$-matrix,
\[
  f_i(\vv_j)=P(\mathbf{M}_{ij})=
  \begin{cases}
    P(0),& \text{if} \ i=j,\\
    0,& \text{if} \ i\ne j.
  \end{cases}
\]
It follows that the polynomials $f_1,\ldots,f_n$ are linearly independent over $K$. Indeed, suppose that $\lambda_1,\dots, \lambda_n\in K$ satisfy
\begin{equation*}
    \sum_{i=1}^n\lambda_if_i = 0.
\end{equation*}
Evaluating at $\vv_j$ gives
$\lambda_jP(0)=0$, so $\lambda_j=0$ for every $1\leq j\leq n$.

The space of polynomials in $r$ variables of total degree at most three has the monomial basis
\[
  \{X_1^{e_1}\cdots X_r^{e_r}:
       e_1,\ldots,e_r\in\mathbb Z_{\geq0},\ e_1+\cdots+e_r\leq3\}.
\]
This basis has size $\binom{r+3}{3}$, so $n\leq\binom{r+3}{3}$.
Finally, $6n\leq (r+1)(r+2)(r+3)\leq (r+3)^3$, which gives $r\geq(6n)^{1/3}-3$.
\end{proof}


\medskip

\begin{thebibliography}{BGGIST07}

\bibitem{bukh-problems}
B.~Bukh.
\emph{Interesting problems that I cannot solve}, Problem~7.
\url{https://www.borisbukh.org/problems.html}.

\bibitem{bukh-ranks}
B.~Bukh.
\emph{Ranks of matrices with few distinct entries}.
Israel J. Math., 222(1):165--200, 2017.

\bibitem{ha}
J.~Ha.
\emph{On the least prime primitive root}.
J. Number Theory, 133:3645--3669, 2013.

\bibitem{lidl}
R.~Lidl and H.~Niederreiter.
\emph{Finite Fields}.
Volume~20 of Encyclopedia of Mathematics and its Applications.
Cambridge University Press, Cambridge, second edition, 1997.

\bibitem{martin}
G.~Martin.
\emph{The least prime primitive root and the shifted sieve}.
Acta Arith., 80:277--288, 1997.


\bibitem{moree}
P.~Moree.
\emph{Artin's primitive root conjecture---a survey}.
Integers, 12(6):1305--1416, 2012.

\bibitem{OpenAI-Artin}
OpenAI.
\href{https://github.com/openai/math/blob/main/preprints/Primitive-roots-for-every-admissible-integer-base-October-4-2026/primitive-roots-all-integer-bases.pdf}
{\emph{Primitive roots for every admissible integer base}}.
OpenAI Math Release preprint, 4 October 2026.

\bibitem{shoup}
V.~Shoup.
\emph{Searching for primitive roots in finite fields}.
Math. Comp., 58:369--380, 1992.

\end{thebibliography}
\end{document}